\documentclass[12pt,leqno,a4paper,dvips,twoside]{article}

\usepackage{graphicx}
\usepackage{amsfonts}
\usepackage{amssymb}
\usepackage{amsmath}
\usepackage{amsthm}
\usepackage{enumerate}                      

\usepackage{calc}                           

\usepackage[utf8]{inputenc}                 

\usepackage{tikz}                           

\newcommand*{\N}{\ensuremath{\mathbb{N}}}               
\newcommand*{\R}{\ensuremath{\mathbb{R}}}               
\newcommand*{\Q}{\ensuremath{\mathbb{Q}}}               

\newtheorem{thm}{Theorem}[section]
\newtheorem{defin}[thm]{Definition}
\newtheorem{lema}[thm]{Lemma}
\newtheorem{cor}[thm]{Corollary}
\newtheorem{conj}[thm]{Conjecture}
\newtheorem{prop}[thm]{Proposition}
\newtheorem{rema}[thm]{Remark}
\newtheorem{exem}[thm]{Example}

\begin{document}
\thispagestyle{plain} {\footnotesize {{}} \vspace*{2cm}

\begin{center}
{\Large {\bf  Chebyshev-Gr\"{u}ss-type inequalities via discrete oscillations}}
\vspace{0.5cm}

{\large Heiner Gonska, Ioan Ra\c{s}a, Maria-Daniela Rusu}\footnote{The authors gratefully acknowledge DAAD support in the framework of the "Stability Pact for Southeastern Europe".}
\end{center}
\vspace{0.5cm}

\begin{abstract}
The classical form of Gr{\"u}ss' inequality, first published by G. Gr\"{u}ss in~\cite{gruess:1935}, gives an estimate of the difference between the integral of the product and the product of the integrals of two functions. In the subsequent years, many variants of this inequality appeared in the literature. The aim of this paper is to introduce a new approach, presenting a new Chebyshev-Gr{\"u}ss-type inequality and applying it to different well-known linear, not necessarily positive, operators. Some interesting conjectures are presented. We also compare the new inequalities with some older results. This new approach gives better estimates in some cases than the ones already known. 

\end{abstract}

\begin{center}

\end{center}
\vspace*{0.3cm}

\section{Introduction}

Here we list some classical results which we will need in the sequel.

The functional given by
\begin{equation*}
T(f,g):=\frac{1}{b-a}\int_{a}^{b}f(x)g(x)dx-\frac{1}{b-a}\int_{a}^{b}f(x)dx\cdot\frac{1}{b-a}\int_{a}^{b}g(x)dx,
\end{equation*}
where $f,g:[a,b]\to \R$ are integrable functions, is well known in the literature as the
classical Chebyshev functional (see \cite{che:fr}).

We first recall the following result. 

\begin{thm}(see~\cite{mitpefi:1993})
	Let $f,g:[a,b]\to\R$ be bounded integrable functions, both increasing or both decreasing. Furthermore, let $p:[a,b]\to\R^{+}_{0}$ be a bounded and integrable function. Then
	\begin{equation}\label{ineq1.1}
		\int^b_a{p(x)dx}\int^b_a{p(x)\cdot f(x)\cdot g(x)dx}\geq\int^b_a{p(x)\cdot f(x)dx}\int^b_a{p(x)\cdot g(x)dx}.
	\end{equation}
If one of the functions $f$ or $g$ is nonincreasing and the other nondecreasing, then inequality ~(\ref{ineq1.1}) is reversed. 
\end{thm}

\begin{rema}
	Inequality ~(\ref{ineq1.1}) is known as Chebyshev's inequality. It was first introduced by P. L. Chebyshev in 1882 in \cite{che:1882}. 
	If $p(x)=1$ for $a\leq x\leq b$, then inequality (\ref{ineq1.1}) is equivalent to 
		\begin{equation*}
			\frac{1}{b-a}\int^{b}_{a}{f(x)\cdot g(x) dx}\geq \left(\frac{1}{b-a}\int^{b}_{a}{f(x) dx}\right)\cdot\left(\frac{1}{b-a}\int^{b}_{a}{g(x) dx}\right).
		\end{equation*}	 
\end{rema}

The next result is the Gr\"{u}ss-type inequality for the Chebyshev functional.

\begin{thm}(Gr\"{u}ss, 1935, see \cite{gruess:1935})\label{t1}
Let $f, g$ be integrable functions from $[a,b]$ into $\R$, such that $m\leq f(x)\leq M$, $p\leq g(x)\leq P$, for
all $x\in [a,b]$, where $m,M,p,P\in\R$. Then
	\begin{equation*}
		\left|T(f,g)\right|\leq \frac{1}{4}(M-m)(P-p).
	\end{equation*}
\end{thm}

The functional $L$, given by $L(f):=\frac{1}{b-a}\int^{b}_{a}{f(x)dx}$, is linear and positive and satisfies $L(e_0)=1$; here we denote $e_i(x)=x^i$, for $i\geq 0$. In the sequel, we recall some bounds for what we call the generalized Chebyshev functional
	\begin{equation}\label{gcf}
		T_L(f,g):=L(f\cdot g)-L(f)\cdot L(g)
	\end{equation}
and give some new results.

\begin{rema}
	We will use the terminology "Chebyshev-Gr\"{u}ss-type inequalities", referring to Gr\"{u}ss-type inequalities for (special cases of) generalized Chebyshev functionals. These inequalities have the general form
		\begin{equation*}	
			\left|T_L(f,g)\right|\leq E(L,f,g),
		\end{equation*}  
	where $E$ is an expression in terms of certain properties of $L$ and some kind of oscillations of $f$ and $g$.
\end{rema} 

Another result we recall is a special form of a theorem given by D. Andrica and C. Badea (see \cite{andrica:1988}):
	\begin{thm}
Let $I=[a,b]$ be a compact interval of the real axis, $B(I)$ be the space of real-valued and bounded
functions defined on $I$ and $L$ be a linear positive functional satisfying $L(e_0)=1$ where $e_0:I\ni x \mapsto 1$. Assuming that for $f,g\in B(I)$ one has $m\leq f(x)\leq M$, $p\leq g(x)\leq P$ for all
$x\in I$, the following holds:
	\begin{equation*}
		\left|T_L(f,g)\right| \leq \frac{1}{4}(M-m)(P-p).
	\end{equation*}
\end{thm}

\begin{rema}
	Note that the positive linear functional is not present on the right hand side of the estimate. 
\end{rema}

The following pre-Chebyshev-Gr\"{u}ss inequality was given by A. Mc. D. Mercer and P. R. Mercer (see \cite{mercer:2004}) in 2004.

\begin{thm}
	For a positive linear functional $L:B(I)\to \R$, with $L(e_0)=1$, one has:
		\begin{equation*}
	\left| T_L(f,g)\right|
		\leq \frac{1}{2}\min\{(M-m) L\left(\left|g-G\right|\right),
			(P-p)L\left(\left|f-F\right|\right)\}
		\end{equation*}	
	where $m\leq f(x)\leq M$, $p\leq g(x)\leq P$ for all $x\in I$, $F:=Lf$ and $G:=Lg$.  
\end{thm}

\begin{rema}
	This is a more adequate result, considering that the positive linear functional appears on both the left and the right hand side of the inequality.
\end{rema}

Let $C(X)=C_{\R}((X,d))$ be the Banach lattice of real-valued continuous functions defined on the compact metric space $(X,d)$ and consider positive linear operators $H:C(X)\to C(X)$ reproducing constant functions. For $x\in X$ we take $L=\epsilon_x\circ H$, so $L(f)=H(f;x)$. We are interested in the degree of non-multiplicativity of such operators. Consider two functions $f,g\in C(X)$ and define
the positive bilinear functional 
	\begin{equation*}
		T(f,g;x):=H(f\cdot g;x)-H(f;x)\cdot H(g;x).
	\end{equation*}

\begin{defin}
	Let $f\in C(X)$. If, for $t\in [0,\infty)$, the quantity
		\begin{equation*}
			\omega_d(f;t):=\sup\left\{\left|f(x)-f(y)\right|,\ d(x,y)\leq t\right\}
		\end{equation*}
	is the usual modulus of continuity, then its least concave majorant is given by
		\begin{align*}
		\widetilde{\omega_d}(f,t)=
			\begin{cases}
\sup_{0\leq x\leq t\leq y\leq d(X),x\neq y}\frac{(t-x)\omega_d(f,y)+(y-t)\omega_d(f,x)}{y-x} & \text{ for } 0\leq t\leq d(X) \text{ , }\\
\omega_d(f,d(X)) & \text{ if } t > d(X)\text{ , }
			\end{cases}
	\end{align*}		
	and $d(X)<\infty$ is the diameter of the compact space $X$.
\end{defin}

In \cite{rusu:2011} (see Theorem 3.1.) the following was shown.
\begin{thm}\label{th2}
	If $f,g\in C(X)$, where $(X,d)$ is a compact metric space, and $x\in X$ fixed, then the inequality
		\begin{equation*}
			\left|T(f,g;x)\right|
			\leq \frac{1}{4}\widetilde{\omega_d}\left(f;4\sqrt{H(d^2(\cdot,x);x)}\right)\cdot
			\widetilde{\omega_d}\left(g;4\sqrt{H(d^2(\cdot,x);x)}\right)
		\end{equation*}	
	holds, where $\widetilde{\omega_d}$ is the least concave majorant of the usual modulus of continuity  		and $H(d^2(\cdot,x);x)$ is the second moment of the operator $H$.
\end{thm}

For $X=[a,b]$, we have a slightly better result (see Theorem 4.1. in \cite{rusu:2011}); a slightly weaker inequality had
been obtained earlier in \cite{acugora:2011}.
\begin{thm}\label{th3}
	If $f,g\in C[a,b]$ and $x\in [a,b]$ is fixed, then the inequality
		\begin{equation}\label{eqth3}
			\left|T(f,g;x)\right|\leq \frac{1}{4}\widetilde{\omega}\left(f;2\sqrt{H((e_1-x)^2;x)}\right)\cdot
			\widetilde{\omega}\left(g;2\sqrt{H((e_1-x)^2;x)}\right)
		\end{equation} 
	holds.
\end{thm}

\begin{rema}
	Here the moduli of continuity are oscillations defined with respect to functions $f$ on the whole domain $X=[a,b]$. In order to improve some results, we propose a new approach, in which the oscillations are related to the support of the involved functional. 
\end{rema}

\begin{rema}
	The inequality (\ref{eqth3}) is sharp in the sense that a positive linear operator reproducing constant and linear functions and functions $f,g\in C[a,b]$ exist such that equality occurs.
\end{rema}

\begin{exem}
	Consider $f=g:=e_1$. Then we have
		\begin{equation*}	
			\omega(f,t)=\omega(e_1,t)=\sup\{\left|x-y\right|:\left|x-y\right|\leq t\}=t. 
		\end{equation*}	
	Since $\omega(f,\cdot)$ is linear, we get $\widetilde{\omega}(f,\cdot)=\omega(f,\cdot)$. 
	The left-hand side in Theorem \ref{th3} is
		\begin{equation*}
			\left|T(f,g;x)\right|=H(e_2;x)-(H(e_1;x))^2
		\end{equation*}
	and the right-hand side is
		\begin{align*}
			\frac{1}{4}\widetilde{\omega}\left(f;2\sqrt{H((e_1-x)^2;x)}\right)\cdot
			\widetilde{\omega}\left(g;2\sqrt{H((e_1-x)^2;x)}\right)&=\frac{1}{4}\cdot (2\sqrt{H((e_1-x)^2;x)})^2\\
			&=H((e_1-x)^2;x).
		\end{align*}
	By choosing a positive linear operator $H:C[a,b]\to [a,b]$ such that $He_0=e_0$ and $He_1=e_1$, we get
		\begin{align*}
			H((e_1-x)^2;x)&=H(e_2-2xe_1+x^2;x)\\
			&=H(e_2;x)-2xH(e_1;x)+x^2=H(e_2;x)-x^2\\
			&=H(e_2;x)-(H(e_1;x))^2,
		\end{align*}
	so we obtain equality between the two sides.
\end{exem}

\section{A Chebyshev-Gr\"{u}ss-type inequality: new approach}

\subsection{The compact topological space case}
Let $\mu$ be a (not necessarily positive) Borel measure on the compact topological space $X$.

Let $\int\limits_X{d\mu(x)}=1$, and set $L(f)=\int\limits_X{f(x)d\mu(x)}$, for $f\in C(X)$. 
Then, for $f,g\in C(X)$, we have

\begin{align*}
	L(fg)-L(f)L(g)&=\int\limits_X{f(x)g(x)d\mu(x)}-\int\limits_X{f(x)d\mu(x)}\cdot \int\limits_X{g(y)d\mu(y)}\\
	&=\iint\limits_{X^2}{f(x)g(x)d(\mu\otimes\mu)(x,y)}-\iint\limits_{X^2}{f(x)g(y)d(\mu\otimes\mu)(x,y)}\\
	&=\iint\limits_{X^2}{f(x)(g(x)-g(y))d(\mu\otimes\mu)(x,y)}.
\end{align*}
Similarly, 
	\begin{equation*}
		L(fg)-L(f) L(g)=\iint\limits_{X^2}{f(y)(g(y)-g(x))d(\mu\otimes\mu)(x,y)}.
	\end{equation*}
By addition,
	\begin{equation}\label{eq5}
		2(L(fg)-L(f) L(g))=\iint\limits_{X^2}{(f(x)-f(y))(g(x)-g(y))d(\mu\otimes\mu)(x,y)}.
	\end{equation}
Let 
	\begin{equation*}
		osc_L(f):=\max\{\left|f(x)-f(y)\right|: (x,y)\in supp (\mu\otimes\mu)\},
	\end{equation*}
where $supp (\mu\otimes\mu)$ is the support of the tensor product of the Borel measure $\mu$ with itself (see \cite{altomare:1994}) and let $\Delta:=\{(x,x):x\in X\}$. From (\ref{eq5}) we get
	\begin{equation*}
		L(f g)-L(f) L(g)=\frac{1}{2}\iint\limits_{X^2\setminus\Delta}{(f(x)-f(y))(g(x)-g(y))d(\mu\otimes\mu)(x,y)}.
	\end{equation*}

Then we have the following result. 

\begin{thm}
	The Chebyshev-Gr\"{u}ss-type inequality in this case is given by
	\begin{equation*}
		\left| L(fg)-L(f) L(g)\right|\leq \frac{1}{2}\cdot osc_L(f)\cdot osc_L(g)\left|\mu\otimes\mu\right|(X^2\setminus \Delta),
	\end{equation*}
	for $f, g\in C(X)$ and $\left|\mu\otimes\mu\right|$ is the absolute value of the tensor product of the Borel measure $\mu$ with itself (see Chapter 1 in \cite{altomare:1994}).
\end{thm}

\begin{exem}
	Let $X=[0,1]$ and consider the functional
	\begin{equation*}
L(f)=a\int^{1}_{0}{f(t)dt}+(1-a)f\left(\frac{1}{2}\right), \text{ for }0\leq a \leq 1.
	\end{equation*}
Then $L(f)=\int^{1}_{0}{f(t)d\mu}$, where the Borel measure $\mu$ is given by 
	\begin{equation*}	
	\mu=a\lambda+(1-a)\varepsilon_{\frac{1}{2}}
	\end{equation*}	
on $X$, with $\lambda$ the Lebesgue measure on $[0,1]$ and $\varepsilon_{\frac{1}{2}}$ the measure concentrated at $\frac{1}{2}$. Then the tensor product of $\mu$ with itself is
\begin{align*}
	\mu\otimes\mu&=\left(a\lambda+(1-a)\varepsilon_{\frac{1}{2}}\right)\otimes \left(a\lambda+(1-a)\varepsilon_{\frac{1}{2}}\right)\\
	&=a^2(\lambda\otimes\lambda)+a(1-a)(\lambda\otimes\varepsilon_{\frac{1}{2}})+(1-a)a(\varepsilon_{\frac{1}{2}}\otimes\lambda)+(1-a)^2(\varepsilon_{\frac{1}{2}}\otimes\varepsilon_{\frac{1}{2}}).
\end{align*}
$\mu\otimes\mu$ is a positive measure, so $\left|\mu\otimes\mu\right|=\mu\otimes\mu$, and 
	\begin{align*}
		&\mu\otimes\mu\left([0,1]^2\setminus \Delta\right)=[a^2(\lambda\otimes\lambda)+a(1-a)(\lambda\otimes \varepsilon_{\frac{1}{2}})\\
		&+a(1-a)(\varepsilon_{\frac{1}{2}}\otimes\lambda)+(1-a)^2 (\varepsilon_{\frac{1}{2}}\otimes \varepsilon_{\frac{1}{2}})]\left([0,1]^2\setminus \Delta\right)\\
		&=a^2+2a(1-a)=a(2-a).
	\end{align*}
The inequality becomes :
	\begin{equation*}
	\left|L(fg)-L(f)L(g)\right|\leq \frac{1}{2}\cdot a(2-a)\cdot osc_L(f)\cdot osc_L(g),
	\end{equation*}
for two functions $f,g\in C[0,1]$.
\end{exem}

\subsection{The discrete linear functional case}

Let $X$ be an arbitrary set and $B(X)$ the set of all real-valued, bounded functions on $X$. Take $a_n\in \R,\ n\geq 0$, such that 
$\sum^{\infty}_{n=0}{\left|a_n\right|}<\infty$ and $\sum^{\infty}_{n=0}{a_n}=1$. Furthermore, let $x_n\in X,\ n\geq 0$ be arbitrary mutually distinct points of $X$. For $f\in B(X)$ set $f_n:=f(x_n)$. Now consider the functional $L:B(X)\to \R$, $Lf=\sum^{\infty}_{n=0}{a_nf_n}$. $L$ is linear and $Le_0=1$. 

Then the relations
	\begin{align*}
		L(f\cdot g)-L(f)\cdot L(g)&=\sum^{\infty}_{n=0}{a_nf_ng_n}-\sum^{\infty}_{n=0}{a_nf_n}\cdot \sum^{\infty}_{m=0}{a_mg_m}\\
		&=\sum^{\infty}_{n=0}{\left(\sum^{\infty}_{m=0}{a_m}\right)a_nf_ng_n}-\sum^{\infty}_{n,m=0}{a_na_mf_ng_m}\\
		&=\sum^{\infty}_{n=0}{a^2_nf_ng_n}+\sum^{\infty}_{n,m=0; m\neq n}{a_ma_nf_ng_n}\\
		&-\sum^{\infty}_{n=0}{a^2_nf_ng_n}-\sum^{\infty}_{n,m=0; m\neq n}{a_na_mf_ng_m}\\
		&=\sum^{\infty}_{n,m=0; m\neq n}{a_na_mf_n(g_n-g_m)}\\
		&=\sum_{0\leq n<m< \infty}{a_na_m f_n(g_n-g_m)}+\sum_{0\leq n>m< \infty}{a_na_mf_n(g_n-g_m)}\\
		&=\sum_{0\leq n<m<\infty}{a_na_mf_n(g_n-g_m)}-\sum_{0\leq n<m< \infty}{a_na_mf_m(g_n-g_m)}\\
		&=\sum_{0\leq n<m< \infty}{a_na_m(f_n-f_m)(g_n-g_m)}
	\end{align*}
hold. 

\begin{thm}\label{th8}
The Chebyshev-Gr\"{u}ss-type inequality for the above linear, not necessarily positive, functional $L$ is given by:
	\begin{equation*}
		\left|L(fg)-L(f)\cdot L(g)\right| \leq osc_L(f)\cdot osc_L(g)\cdot \sum_{0\leq n<m< \infty}{\left|a_na_m\right|},
	\end{equation*}
where $f,g\in B(X)$  and we define the oscillations to be:
	\begin{align*}
		osc_L(f)&:=\sup\{\left|f_n-f_m\right|: 0\leq n<m< \infty\},\\
		osc_L(g)&:=\sup\{\left|g_n-g_m\right|: 0\leq n<m< \infty\}.
	\end{align*}
\end{thm}

\begin{thm}\label{th9}
In particular, if $a_n\geq 0$, $n\geq 0$, then $L$ is a positive linear functional and we have:
	\begin{equation*}
		\left|L(fg)-Lf\cdot Lg\right| \leq \frac{1}{2}\cdot \left(1-\sum^{\infty}_{n=0}{a^2_n}\right)\cdot osc_L(f)\cdot osc_L(g),
	\end{equation*}
for $f,g\in B(X)$ and the oscillations given as above.
\end{thm}

\begin{rema}
	The above inequality is sharp in the sense that we can find a functional $L$ such that equality holds. 
\end{rema}

\begin{exem}
	Let us consider the following functional 
		\begin{equation*}
			Lf:=(1-a)f(0)+af(1), \text{ for } 0\leq a\leq 1.
		\end{equation*}
	For this functional we have
		\begin{equation*}
			L(fg)-Lf\cdot Lg=(1-a)f(0)g(0)+af(1)g(1)-[(1-a)f(0)+af(1)]\cdot[(1-a)g(0)+ag(1)],
		\end{equation*}
	so after some calculations we get that the left-hand side is
		\begin{align*}
			\left|L(fg)-Lf\cdot Lg\right|&=\left|\underbrace{a(1-a)}_{\geq 0}\cdot[f(0)-f(1)]\cdot[g(0)-g(1)]\right|\\
			&=a(1-a)\left|f(0)-f(1)\right|\cdot\left|g(0)-g(1)\right|
		\end{align*}
	and the right-hand side is
	\begin{align*}
			\frac{1}{2}\left(1-\sum^{\infty}_{n=0}{a^2_n}\right)\cdot osc_L(f)\cdot osc_L(g)&=\frac{1}{2}\cdot [1-a^2-(1-a)^2]\cdot \left|f(0)-f(1)\right|\cdot\left|g(0)-g(1)\right|\\
			&=a(1-a)\left|f(0)-f(1)\right|\cdot\left|g(0)-g(1)\right|.
	\end{align*}
\end{exem}

\section{A new Chebyshev-Gr\"{u}ss-type inequality for the Bernstein operator}

Consider the classical Bernstein operators
	\begin{equation*}
		B_nf(x):=\sum^{n}_{k=0}{f\left(\frac{k}{n}\right)b_{nk}(x)},\ f\in \R^{[0,1]},\ x\in [0,1],
	\end{equation*}
	where $b_{nk}(x):=\binom{n}{k}x^k(1-x)^{n-k}$.
According to Theorem \ref{th9}, for each $x\in [0,1]$, $f,g\in B[0,1]$ we have
	\begin{equation}\label{eqalpha}
		\left|B_n(f\cdot g)(x)-B_nf(x)\cdot B_ng(x)\right|\leq\frac{1}{2}\left(1-\sum^{n}_{k=0}{b^2_{nk}(x)}\right)\cdot osc_{B_n}(f)\cdot osc_{B_n}(g),
	\end{equation}	
where 
	\begin{equation*}
		osc_{B_n}(f):=\max\{\left|f_k-f_l\right|:0\leq k<l\leq n\}
	\end{equation*}	 
and $f_k:=f\left(\frac{k}{n}\right)$; similar definitions apply to $g$. 

\begin{exem}
If we consider $f,g\in B[0,1]$ to be Dirichlet functions defined by
		\begin{equation*}
			f(x):=\begin{cases}
				1 & \text{ for } x\in \Q,\\
				0 & \text{ for }x\in \R\setminus\Q
			\end{cases}
		\end{equation*}
and analogously for $g$,  with $f_k:=f\left(\frac{k}{n}\right)$ (the same for $g$), then we observe that 	the oscillations in the above inequality vanish, so the right hand-side is zero. 	
\end{exem}

Let $\varphi_n(x):=\sum^{n}_{k=0}{b^2_{nk}(x)}$, $x\in [0,1]$. Since
	\begin{equation*}
		\left(\frac{1}{n+1}\sum^{n}_{k=0}{b^2_{nk}(x)}\right)^{\frac{1}{2}}\geq \frac{1}{n+1}\sum^{n}_{k=0}{b_{nk}(x)}=\frac{1}{n+1},
	\end{equation*}	
we get
	\begin{equation}\label{eq2}
		\varphi_n(x)\geq \frac{1}{n+1},\ x\in [0,1],
	\end{equation}
and therefore
	\begin{equation}\label{eqbeta}
		\left|B_n(f\cdot g)(x)-B_nf(x)\cdot B_ng(x)\right|\leq\frac{n}{2(n+1)}\cdot osc_{B_n}(f)\cdot osc_{B_n}(g),\ x\in [0,1].
	\end{equation}	
Let us remark that equality is attained in (\ref{eq2}) iff $n=1$ and $x=\frac{1}{2}$. In fact, inspired also by some computations with Maple, we
make the following conjectures:

\begin{conj}\label{c1}
	$\varphi_n$ is convex on $[0,1]$.
\end{conj}

\begin{conj}\label{c2}
	$\varphi_n$ is decreasing on $\left[0,\frac{1}{2}\right]$ and increasing on $\left[\frac{1}{2},1\right]$.
\end{conj}

\begin{conj}\label{c3}
	$\varphi_n(x)\geq \varphi_n\left(\frac{1}{2}\right)$, $x\in [0,1]$.
\end{conj}

Since $\varphi_n\left(\frac{1}{2}-t\right)=\varphi_n\left(\frac{1}{2}+t\right)$, $t\in \left[0,\frac{1}{2}\right]$, we see that 

Conjecture \ref{c1} $\Rightarrow$ Conjecture \ref{c2} $\Rightarrow$ Conjecture \ref{c3}.

On the other hand, it can be proved that
	\begin{equation*}
		\varphi_n\left(\frac{1}{2}\right)=4^{-n}\binom{2n}{n},\ \varphi'_n\left(\frac{1}{2}\right)=0,\ \varphi''_n\left(\frac{1}{2}\right)=4^{2-n}\binom{2n-2}{n-1},
	\end{equation*}
and so $\frac{1}{2}$ is a minimum point for $\varphi_n$. Conjecture \ref{c3} claims that it is an absolute minimum point; in other words,
	\begin{equation}\label{eq3}
		\varphi_n(x)\geq \frac{1}{4^n}\binom{2n}{n}, \ x\in [0,1].
	\end{equation}

The following confirmation of Conjecture \ref{c3} is due to Dr. Th. Neuschel (University of Trier). 

\begin{lema}\label{lc3}
	For $n\in\N$ and $x\in [0,1]$, we have
		\begin{equation*}
			\sum^{n}_{k=0}{\binom{n}{k}}^2x^{2k}(1-x)^{2(n-k)}\geq \frac{1}{4^n}\binom{2n}{n}.
		\end{equation*}
\end{lema}

\begin{proof}
	For symmetry reasons, it suffices to prove the statement only for $0\leq x\leq\frac{1}{2}$. 
	In the sequel we denote $P_n$ to be the $n-$th Legendre polynomial, given by
		\begin{equation*}
			P_n(x):=\frac{1}{2^n}\sum^{n}_{k=0}{\binom{n}{k}}^2(x+1)^k(x-1)^{n-k}.
		\end{equation*}
	We make a change of variable, namely set $y:=\frac{1-2x+2x^2}{1-2x}\geq 1$ and we get
		\begin{equation*}
			(y-\sqrt{y^2-1})^n\cdot P_n(y)=\sum^{n}_{k=0}{\binom{n}{k}}^2x^{2k}(1-x)^{2(n-k)}:=\varphi_n(x),
		\end{equation*}
	so we have to show that
		\begin{equation*}
			(y-\sqrt{y^2-1})^n\cdot P_n(y)\geq \frac{1}{4^n}\binom{2n}{n}
		\end{equation*}
	holds, for $y\geq 1$. The inequality holds for $y=1$ and $y=\infty$. In the last case, the inequality is even sharp. 
	Now it is enough to show :
		\begin{equation*}
			\frac{d}{dy}\{(y-\sqrt{y^2-1})^nP_n(y)\}\leq 0 \text{ for } y>1. 
		\end{equation*}
	This is equivalent to the following statement
		\begin{equation*}
			P'_n(y)\leq \frac{n}{\sqrt{y^2-1}}P_n(y) \text{ for } y>1.
		\end{equation*}
	Using the formula
		\begin{equation*}
			\frac{y^2-1}{n}P'_n(y)=yP_n(y)-P_{n-1}(y),
		\end{equation*}
	we now have to prove the following:
		\begin{equation*}
			(y-\sqrt{y^2-1})P_n(y)\leq P_{n-1}(y) \text{ for } y>1,
		\end{equation*}
	which is equivalent to
		\begin{equation}\label{eq0.3}
			P_n(y)\leq (y+\sqrt{y^2-1})P_{n-1}(y) \text{ for } y>1.
		\end{equation}
	The inequality (\ref{eq0.3}) can be proved by induction. For $n=1$ the inequality holds. We assume that the inequality holds also for $n$ and 
	we want to show:
		\begin{equation*}
			P_{n+1}(y)\leq (y+\sqrt{y^2-1})P_n(y) \text{ for } y>1.
		\end{equation*}
	Using Bonnet's recursion formula
		\begin{equation*}
			P_{n+1}(y)=\frac{2n+1}{n+1}yP_n(y)-\frac{n}{n+1}P_{n-1}(y),
		\end{equation*}
	we now have to show that the following holds:
		\begin{equation*}
			\left(\frac{2n+1}{n+1}y-(y+\sqrt{y^2-1})\right)P_n(y)\leq \frac{n}{n+1}P_{n-1}(y).
		\end{equation*}
	After evaluation
		\begin{align*}
			\left(\frac{2n+1}{n+1}y-(y+\sqrt{y^2-1})\right)P_n(y)&\leq \frac{n}{n+1}(y-\sqrt{y^2-1})P_n(y)\\
			&\leq \frac{n}{n+1}(y-\sqrt{y^2-1})(y+\sqrt{y^2-1})P_{n-1}(y)\\
			&=\frac{n}{n+1}P_{n-1}(y),
		\end{align*}
	we obtain the result. 
\end{proof}

In order to compare (\ref{eq2}) and (\ref{eq3}), it is not difficult to prove the inequalities
	\begin{equation*}
		\frac{1}{n+1}< \frac{1}{2\sqrt{n}}< \frac{1}{4^n}\binom{2n}{n}< \frac{1}{\sqrt{2n+1}},\ n\geq 2.
	\end{equation*}
More precise inequalities can be found in \cite{elpre:2012}:
	\begin{equation*}
		\frac{1}{\sqrt{\pi(n+3)}}< \frac{1}{4^n}\binom {2n}{n}< \frac{1}{\sqrt{\pi(n-1)}}, \ n\geq 2. 
	\end{equation*}
	
Because we have proved that Conjecture \ref{c3} is true, we have the following result.
\begin{thm} 
The new Chebyshev-Gr\"{u}ss-type inequality for the Bernstein operator is:
	\begin{equation}\label{eqgamma}
		\left|B_n(f\cdot g)(x)-B_nf(x)\cdot B_ng(x)\right|\leq\frac{1}{2}\left(1-\frac{1}{4^n}\binom{2n}{n}\right)\cdot osc_{B_n}(f)\cdot osc_{B_n}(g),\ x\in[0,1].
	\end{equation}
\end{thm}

In comparison, using the second moment of the Bernstein polynomial
	\begin{equation*}
		B_n((e_1-x)^2;x)=\frac{x(1-x)}{n},
	\end{equation*}
and letting $H=B_n$ in Theorem \ref{th3}, the classical Chebyshev-Gr\"{u}ss-type inequality looks as follows:
	\begin{equation}\label{eqeta}
		\left|B_n(f\cdot g)(x)-B_nf(x)\cdot B_ng(x)\right|\leq \frac{1}{4}\widetilde{\omega}\left(f;2\sqrt{\frac{x(1-x)}{n}}\right)\cdot\widetilde{\omega}\left(g;2\sqrt{\frac{x(1-x)}{n}}\right),
	\end{equation}
	which implies
	\begin{equation}\label{eqepsilon}
	\left|B_n(f\cdot g)(x)-B_nf(x)\cdot B_ng(x)\right|\leq \frac{1}{4}\widetilde{\omega}\left(f;\frac{1}{\sqrt{n}}\right)\cdot\widetilde{\omega}\left(g;\frac{1}{\sqrt{n}}\right),	
	\end{equation} 
for two functions $f,g\in C[0,1]$ and $x\in [0,1]$ fixed. 

\begin{rema}
	In (\ref{eqalpha}) and (\ref{eqeta}), the right-hand side depends on $x$ and vanishes when $x\to 0$ or $x\to 1$. The maximum
	value of it, as a function of $x$, is attained for $x=\frac{1}{2}$, and (\ref{eqbeta}), (\ref{eqgamma}), (\ref{eqepsilon}) illustrate this fact.
	On the other hand, in (\ref{eqalpha}) the oscillations of $f$ and $g$ are relative only to the points $0,\frac{1}{n},\ldots, \frac{n-1}{n},1$, while in
	(\ref{eqeta}) the oscillations, expressed in terms of $\widetilde{\omega}$, are relative to the whole interval $[0,1]$. 
\end{rema}

\section{Gr\"{u}ss-type inequalities for the Lagrange operator}

Consider $f\in C[-1,1]$ and the infinite matrix $X=\{x_{k,n}\}^{n\ \ \ \ \infty}_{k=1 \ n=1}$ with
	\begin{equation*}
		-1\leq x_{1,n}<x_{2,n}<\ldots <x_{n,n}\leq 1, \text{ for } n=1,2,\ldots.
	\end{equation*}
The Lagrange fundamental functions are given by
	\begin{equation*}
		l_{k,n}(x)=\frac{\omega_n(x)}{\omega'_n(x_{k,n})(x-x_{k,n})},\ 1\leq k\leq n,
	\end{equation*}
where $\omega_n(x)=\prod^{n}_{k=1}{(x-x_{k,n})}$ and the Lagrange operator (see \cite{szaver:1990}) $L_n:C[-1,1]\to \Pi_{n-1}$ is
	\begin{equation*}
		L_n(f;x):=\sum^{n}_{k=1}{f(x_{k,n})l_{k,n}(x)}. 
	\end{equation*}
The Lebesgue function of the interpolation is:
	\begin{equation*}
		\Lambda_n(x):=\sum^{n}_{k=1}{\left|l_{k,n}(x)\right|}.
	\end{equation*}
It is also known (see \cite{chelight:2000}, p. 13) that $\left\|L_n\right\|<\infty$ and 
	\begin{equation*}
		\left\|L_n\right\|=\left\|\Lambda_n\right\|_{\infty}
	\end{equation*}
hold. 

\begin{prop}[Properties of the Lagrange operator]\mbox{}

	\begin{enumerate}[i)]
		\item
			The Lagrange operator is linear but only in exceptional cases positive.
		\item
			$L_n( f;x_{k,n})=f(x_{k,n}),\ 1\leq k\leq n$.
		\item
			The Lagrange operator is idempotent: $L^2_n=L_n$.
		\item
			$L_n$ satisfies $\sum^{n}_{k=1}{l_{k,n}(x)}=1$. 
	\end{enumerate}
\end{prop}

\begin{rema}
	The Lebesgue function has been studied for different node systems. In the sequel, we will use
	some known results for Chebyshev nodes and give classical and new Chebyshev-Gr\"{u}ss-type inequalities. 
\end{rema}

\subsection{A Chebyshev-Gr\"{u}ss-type inequality for the Lagrange operator at Chebyshev nodes}

The Lagrange operator with Chebyshev nodes (see \cite{brutman:1997}, \cite{chelight:2000}) is given as follows. 

Let $T_n(x)=\cos(n\cos^{-1}x)$ and $X=\{\cos[\pi(2k-1)/2n]\}$, i.e., when 
	\begin{equation*}
		x_{k,n}=\cos t_{k,n}=\cos \frac{2k-1}{2n}\cdot \pi\ (k=1,2,\ldots,n; n=1,2,\ldots)
	\end{equation*} 
are the Chebyshev roots. 

\begin{rema}
It can be shown that the Lebesgue constant for Chebyshev nodes is a lot smaller than for equidistant nodes. That's why we concentrate on this case in our paper. 
\end{rema}

A Chebyshev-Gr\"{u}ss-type inequality for the Lagrange operator with this node system, similar to the one in Theorem \ref{th2}, is given by:

\begin{thm}
	For $f,g\in C[-1,1]$ and all $x\in [-1,1]$, the inequality
		\begin{align*}
			\left| T(f,g;x)\right|
&\leq \frac{1}{4}\left\|L_n\right\|(1+\left\|L_n\right\|)\widetilde{\omega}\left(f;2\right)\cdot \widetilde{\omega}\left(g;2\right)\\
&\leq \frac{1}{2}\left(1+\frac{3}{\pi}\log n+\frac{2}{\pi}\log^2 n\right)\omega(f;2)\cdot \omega(g;2)
		\end{align*}
holds; here $\omega$ denotes the first order modulus. 
\end{thm}

\begin{proof}
The idea of this proof is similar to the one of Theorem 2 in \cite{acugora:2011} and that of Theorem 3.1. in \cite{rusu:2011}. Recall, however, 
that we have to work without the assumption of positivity.
We consider the bilinear functional 
	\begin{equation*}
		T(f,g;x):=L_n(f\cdot g;x)-L_n(f;x)\cdot L_n(g;x).
	\end{equation*}
Let $f,g\in C[-1,1]$ and $r,s\in Lip_1$, where $Lip_1=\{f\in C[-1,1]: \sup_{x\neq x_0}\frac{\left|f(x)-f(x_0)\right|}{\left|x-x_0\right|}<\infty\}$ and the seminorm on $Lip_1$ is defined by $\left|f\right|_{Lip_1}:=\sup_{x\neq x_0}\frac{\left|f(x)-f(x_0)\right|}{\left|x-x_0\right|}$. 
We are interested in estimating
\begin{align} \label{ineq*}
	\begin{split}
		\left|T(f,g;x)\right| &=\left|T(f-r+r,g-s+s;x)\right|\\
								       &\leq \left|T(f-r,g-s;x)\right|+\left|T(f-r,s;x)\right|+\left|T(r,g-s;x)\right|+\left|T(r,s;x)\right|.	
	\end{split}	
\end{align}
First note that for $f,g\in C[-1,1]$ one has
	\begin{equation*}
		\left|T(f,g;x)\right|\leq \left\|L_n\right\|(1+\left\|L_n\right\|)\left\|f\right\|\cdot \left\|g\right\|.
	\end{equation*}
For $r,s\in Lip_1$ we have the estimate
	\begin{align*}
		\left|T(r,s;x)\right|&=\left|T((r-r(0)),(s-s(0);x)\right|\\
		&=\left|L_n((r-r(0))\cdot (s-s(0));x)-L_n(r-r(0);x)\cdot L_n(s-s(0);x)\right|\\
		&\leq \left\|L_n\right\|\cdot \left\|r-r(0)\right\|\cdot \left\|s-s(0)\right\|+\left\|L_n\right\|^2\cdot \left\|r-r(0)\right\|\cdot  \left\|s-s(0)\right\|\\
		&\leq \left\|L_n\right\|(1+\left\|L_n\right\|)\cdot\left|r\right|_{Lip_1}\cdot\left|s\right|_{Lip_1}.
	\end{align*}

Moreover, for $r\in Lip_1$ and $g\in C[-1,1]$ the inequality
	\begin{align*}
		\left|T(r,g;x)\right|&=\left|T(r-r(0),g;x)\right|\\
								   &=\left|L_n((r-r(0))\cdot g;x)-L_n(r-r(0);x)\cdot L_n(g;x)\right|\\
								   &\leq \left\|L_n\right\|\cdot\left\|(r-r(0))\cdot g\right\|+\left\|L_n\right\|^2\cdot \left\|r-r(0)\right\|\cdot\left\|g\right\|\\
								   &\leq \left\|L_n\right\|(1+\left\|L_n\right\|)\cdot\left\|g\right\|\cdot \left\|r-r(0)\right\|\\
								   &\leq \left\|L_n\right\|(1+\left\|L_n\right\|)\cdot\left\|g\right\|\cdot \left|r\right|_{Lip_1}
	\end{align*}
holds. Note that in both cases considered so far we used
	\begin{align*}
			\left|r(x)-r(0)\right|&=\frac{\left|r(x)-r(0)\right|}{\left|x-0\right|}\cdot \left|x-0\right|\\
			&\leq \left|r\right|_{Lip_1}\cdot \left|x\right|,	
	\end{align*} 
for $x\in [-1,1]$, i.e.,  
	\begin{equation*}
		\left\|r(x)-r(0)\right\|\leq \left|r\right|_{Lip_1}.
	\end{equation*}
Similarly, if $f\in C[-1,1]$ and $s\in Lip_1$ we have
	\begin{equation*}
		\left|T(f,s;x)\right|\leq \left\|L_n\right\|(1+\left\|L_n\right\|)\cdot \left\|f\right\|\cdot \left|s\right|_{Lip_1}.
	\end{equation*}
Then inequality (\ref{ineq*}) becomes
	\begin{align*}
		\left|T(f,g;x)\right|
		&\leq \left|T(f-r,g-s;x)\right|+\left|T(f-r,s;x)\right|+\left|T(r,g-s;x)\right|+\left|T(r,s;x)\right|\\
		&\leq \left\|L_n\right\|(1+\left\|L_n\right\|)\cdot\left\{\left\|f-r\right\|+\left|r\right|_{Lip_1}\right\}\cdot\left\{\left\|g-s\right\|+\left|s\right|_{Lip_1}\right\}.
	\end{align*} 
The latter expression involves terms figuring in the K - functional
	\begin{align*}
		&K(f,t;C[-1,1],Lip_1)\\
		&=\inf\{\left\|f-g\right\|+t\cdot \left|g\right|_{Lip_1}: g\in Lip_1\},	
	\end{align*} 
for $f\in C[-1,1]$, $t\geq 0$.
It is known that (see, e.g., \cite{paltanea:2010})
	\begin{equation*}
		K\left(f,\frac{t}{2}\right)=\frac{1}{2}\cdot \widetilde{\omega}(f;t),
	\end{equation*}
an equality to be used in the next step.

We now pass to the infimum over $r$ and $s$, respectively, which leads to
	\begin{align*}
		\left|T(f,g;x)\right|
		&\leq \left\|L_n\right\|(1+\left\|L_n\right\|)\cdot K(f,1;C,Lip_1)\cdot K(g,1;C,Lip_1)\\
		&=\left\|L_n\right\|(1+\left\|L_n\right\|)\cdot\frac{1}{2}\cdot \widetilde{\omega}(f;2)\cdot\frac{1}{2}\cdot\widetilde{\omega}(g;2)\\
		&=\frac{1}{4}\left\|L_n\right\|(1+\left\|L_n\right\|)\omega(f;2)\cdot\omega(g;2).
	\end{align*}
	T. Rivlin (see \cite{rivlin:1981}) proved the following inequality in the case of Lagrange interpolation at Chebyshev nodes:	
	\begin{equation*}
		0.9625<\left\|L_n\right\|-\frac{2}{\pi}\log n<1,
	\end{equation*}
  so using this result we get
		\begin{align*}
			&\left\|L_n\right\|<\frac{2}{\pi}\log n+1\\
			&\Rightarrow 1+\left\|L_n\right\|<2\left(\frac{1}{\pi}\log n+1\right)\\
			&\Rightarrow \left\|L_n\right\|(1+\left\|L_n\right\|)< 2\left(1+\frac{3}{\pi}\log n+\frac{2}{\pi^2}\log^2 n\right)
		\end{align*}
	which implies the result. 
\end{proof}

\subsection{A new Chebyshev-Gr\"{u}ss-type inequality for the Lagrange operator with Chebyshev nodes}

\begin{thm}
For $f,g\in C[-1,1]$ and $x\in [-1,1]$ fixed, the following inequality
		\begin{align*}
	\left| T(f,g;x)\right|
	&\leq osc_{L_n}(f)\cdot osc_{L_n}(g)\cdot \sum_{1\leq k<m\leq n}{\left|l_{k,n}(x)\cdot l_{m,n}(x)\right|}\\
&\leq osc_{L_n}(f)\cdot osc_{L_n}(g)\cdot \left\{\frac{\Lambda^2_n(x)-c\left[1+\cos^{2}nt\cdot\frac{\pi^2}{6}\right]}{2}\right\}
		\end{align*}
holds, for a suitable constant $c$ and $x=\cos t$. 
\end{thm}

\begin{proof}
	The first inequality follows from Theorem \ref{th8} (with an obvious modification). 
	The sum on the right-hand side of the first inequality can  be expressed as follows:
	\begin{align*}
	\sum_{1\leq k<m\leq n}{\left|l_{k,n}(x)\cdot l_{m,n}(x)\right|}&=\left[\left(\sum^{n}_{i=1}{\left|l_{i,n}(x)\right|}\right)^2-\left(\sum^{n}_{i=1}{l^2_{i,n}(x)}\right)\right]/2\\
	&=\left[\Lambda^2_n(x)-\left(\sum^{n}_{i=1}{l^2_{i,n}(x)}\right)\right]/2.
	\end{align*}
In order to estimate the sum $\sum^{n}_{i=1}{l^2_{i,n}(x)}$, we use the proof of Theorem 2.3. from \cite{hermann:1981} to get (case $\alpha=2$):
	\begin{equation*}
		\sum^{n}_{i=1}{l^2_{i,n}(x)}\geq c\left(1+\left|\cos nt\right|^2\sum ^{n}_{i=1}{i^{-2}}\right),  
	\end{equation*}
where $x=\cos t$ and $c$ is a suitable constant. 
After some calculation, the sum becomes
\begin{equation*}
	\sum_{1\leq k<m\leq n}{\left|l_{k,n}(x)\cdot l_{m,n}(x)\right|}
	=\frac{\Lambda^2_n(x)}{2}-\frac{c\left(1+(\cos nt)^2\cdot \frac{\pi ^2}{6}\right)}{2},
\end{equation*}
so we obtain our desired inequality.	
\end{proof}

\section{Chebyshev-Gr\"{u}ss-type inequalities for piecewise linear interpolation at equidistant knots}

We consider the operator $S_{\Delta_n}:C[0,1]\to C[0,1]$ (see \cite{gonska:2011}) at the points $0,\frac{1}{n},\ldots, \frac{k}{n},\ldots, \frac{n-1}{n},1$, which can be explicitely described as 
	\begin{equation*}
		S_{\Delta_n}(f;x)=\frac{1}{n}\sum^{n}_{k=0}{\left[\frac{k-1}{n},\frac{k}{n},\frac{k+1}{n};\left|\alpha-x\right|\right]_{\alpha} f\left(\frac{k}{n}\right)},
	\end{equation*}
where $[a,b,c;f]=[a,b,c;f(\alpha)]_{\alpha}$ denotes the divided difference of a function $f:D\to\R$ on the
(distinct knots) $\{a,b,c\}\subset D$, w.r.t. $\alpha$.

\begin{prop}[Properties of $S_{\Delta_n}$]
	\begin{enumerate}[i)]\mbox{}
		\item
			$S_{\Delta_n}$ is a positive, linear operator preserving linear functions.
		\item
			$S_{\Delta_n}$ preserves monotonicity and convexity/concavity. 
		\item
			$S_{\Delta_n}(f;0)=0,\ S_{\Delta_n}(f;1)=f(1)$.
		\item
			If $f\in C[0,1]$ is convex, then $S_{\Delta_n}f$ is also convex and we have: $f\leq S_{\Delta_n}f$.
	\end{enumerate}
\end{prop}

The operator $S_{\Delta_n}$ can also be defined as follows. 
	\begin{equation*}
		S_{\Delta_n}f(x):=\sum^{n}_{k=0}{f\left(\frac{k}{n}\right)u_{n,k}(x)},
	\end{equation*}
for $f\in C[0,1]$ and $x\in [0,1]$, where $u_{n,k}\in C[0,1]$ are piecewise linear continuous functions, such that
	\begin{equation*}
		u_{n,k}\left(\frac{l}{n}\right)=\delta_{kl},\ k,l=0,\ldots, n. 
	\end{equation*} 

\subsection{A Chebyshev-Gr\"{u}ss-type inequality for $S_{\Delta_n}$}
	
In order to obtain a classical Chebyshev-Gr\"{u}ss-type inequality using $S_{\Delta_n}$, we need the second moment of the operator.
For $x\in \left[\frac{k-1}{n},\frac{k}{n}\right]$, this is given by
	\begin{align*}
		S_{\Delta_n}((e_1-x)^2;x)&=n\left(x-\frac{k-1}{n}\right)\left(\frac{k}{n}-x\right)\left[\left(\frac{k}{n}-x\right)-\left(\frac{k-1}{n}-x\right)\right]\\
		&=\left(x-\frac{k-1}{n}\right)\left(\frac{k}{n}-x\right),
	\end{align*}
which is maximal when $x=\frac{2k-1}{2n}$. This implies
	\begin{equation*}
		S_{\Delta_n}((e_1-x)^2;x)\leq \frac{1}{4n^2}.
	\end{equation*}

By taking $H=S_{\Delta_n}$ in Theorem \ref{th3}, the Chebyshev-Gr\"{u}ss-type inequality for $S_{\Delta_n}$ is given in the following.

\begin{thm}
	If $f,g\in C[0,1]$ and $x\in [0,1]$ is fixed, then the inequality
	\begin{align*}
		\left|T(f,g;x)\right|&\leq \frac{1}{4}\widetilde{\omega}\left(f;2\cdot\sqrt{S_{\Delta_n}((e_1-x)^2;x)}\right)\cdot \widetilde{\omega}\left(g;2\cdot\sqrt{S_{\Delta_n}((e_1-x)^2;x)}\right)\\
		&\leq \frac{1}{4}\widetilde{\omega}\left(f;\frac{1}{n}\right)\cdot\widetilde{\omega}\left(g;\frac{1}{n}\right)
	\end{align*}
	holds. 
\end{thm}

\subsection{A new Chebyshev-Gr\"{u}ss-type inequality for $S_{\Delta_n}$}

In this case, we need to find the minimum of the sum $\tau_n(x):=\sum^{n}_{k=0}{u^2_{n,k}}$. For a particular interval $\left[\frac{k-1}{n},\frac{k}{n}\right]$, we get that 
	\begin{align*}
		\tau_n(x)&:=\sum^{n}_{k=0}{u^2_{n,k}}\\
				  &=(nx-k+1)^2+(k-nx)^2, \text{ for } k=1,\ldots, n. 
	\end{align*}
For $k=1$, we have $x\in \left[0,\frac{1}{n}\right]$ and $\tau_n(x)=(nx-1)^2$, while for $k=n$, we get $x\in \left[\frac{n-1}{n},1\right]$ and $\tau_n(x)=(nx-n+1)^2$. So $\tau_n(x)=(nx-k+1)^2+(k-nx)^2$ is minimal if and only if
$x=\frac{2k-1}{2n}$ and the minimum value of $\tau_n(x)$ is $\frac{1}{2}$. 

\begin{thm}
The new Chebyshev-Gr\"{u}ss-type inequality for $S_{\Delta_n}$ is 
	\begin{align*}
		\left|T(f,g;x)\right|
		&\leq\frac{1}{2}\left(1-\sum^{n}_{k=0}{u^2_{n,k}(x)}\right)\cdot osc_{S_{\Delta_n}}(f)\cdot osc_{S_{\Delta_n}}(g)\\
		&\leq \frac{1}{2}\left(1-\frac{1}{2}\right)\cdot osc_{S_{\Delta_n}}(f)\cdot osc_{S_{\Delta_n}}(g)\\
		&\leq \frac{1}{4}\cdot osc_{S_{\Delta_n}}(f)\cdot osc_{S_{\Delta_n}}(g),
	\end{align*}

with 	
	\begin{align*}
		osc_{S_{\Delta_n}}(f)&:=\max \left\{\left|f_{k}-f_{l}\right|:0\leq k<l\leq n\right\},\\
		osc_{S_{\Delta_n}}(g)&:=\max \left\{\left|g_{k}-g_{l}\right|:0\leq k<l\leq n\right\},
	\end{align*}
where $f_k:=f\left(\frac{k}{n}\right)$.
\end{thm}

\begin{rema}
	This inequality implies the classical Chebyshev-Gr\"{u}ss-type inequality because
	$\left|f_k-f_l\right|\leq M-m$ and $\left|g_k-g_l\right|\leq P-p$, respectively.  It is easy to give examples in which our
	approach gives strictly better inequalities. 	
\end{rema}

\section{Chebyshev-Gr\"{u}ss-type inequalities for \\ Mirakjan-Favard-Sz$\acute{a}$sz operators}

The Mirakjan-Favard-Sz$\acute{a}$sz operators (see \cite{altomare:1994}) were introduced by G. M. Mirakjan (see \cite{mir:1941}) and studied by different authors, e.g.,  J. Favard and O. Sz$\acute{a}$sz (see \cite{favard:1944} and \cite{szasz:1950}). The classical $n-$th Mirakjan-Favard-Sz$\acute{a}$sz operator $M_n$ is defined by
	\begin{equation}\label{smeq1}
		M_n(f;x):=e^{-nx}\sum^{\infty}_{k=0}{\frac{(nx)^k}{k!}f\left(\frac{k}{n}\right)},
	\end{equation}
for $f\in E_2$, $x\in [0,\infty)\subset\R$ and $n\in \N$. $E_2$ is the Banach lattice 
	\begin{equation*}
		E_2:=\{f\in C([0,\infty)): \frac{f(x)}{1+x^2} \text{ is convergent as }x\to\infty\},
	\end{equation*}
endowed with the norm
	\begin{equation*}
		\left\|f\right\|_{*}:=\sup_{x\geq 0}\frac{\left|f(x)\right|}{1+x^2}.
	\end{equation*}
The series on the right-hand side of (\ref{smeq1}) is absolutely convergent and $E_2$ is isomorphic to $C[0,1]$; (see \cite{{altomare:1994}}, Sect. 5.3.9).  

\subsection{A new Chebyshev-Gr\"{u}ss-type inequality for Mirakjan-Favard-Sz$\acute{a}$sz operators}
This is our first application of Theorem \ref{th9} for operators defined for functions given on an infinite interval. We set
	\begin{equation*}
		\sigma_n(x):=e^{-2nx}\sum^{\infty}_{k=0}{\frac{(nx)^{2k}}{(k!)^2}}
	\end{equation*}
and we want to find the infimum:
	\begin{equation*}	
	\inf_{x\geq 0}\sigma_n(x):=\iota\geq 0.
	\end{equation*}
Because $\sigma_n(x)\geq \iota$, we obtain the following result.

\begin{thm}
	For the Mirakjan-Favard-Sz$\acute{a}$sz operator we have
	\begin{align*}
		\left|T(f,g;x)\right|&\leq \frac{1}{2}\left(1-\sigma_n(x)\right)\cdot osc_{M_n}(f)\cdot osc_{M_n}(g)\\
		&\leq \frac{1}{2}(1-\iota)\cdot osc_{M_n}(f)\cdot osc_{M_n}(g),
	\end{align*}	
where $f,g\in C_b[0,\infty)$, $osc_{M_n}(f)=\sup\{\left|f_k-f_l\right|:0\leq k<l<\infty\}$, with $f_k:=f\left(\frac{k}{n}\right)$ and a similar definition applying to $g$. $C_b[0,\infty)$ is the set of all continuous, real-valued, bounded functions on $[0,\infty)$.	
\end{thm}

\begin{lema}
	The relation
	\begin{equation*}
		\inf_{x\geq 0}\sigma_n(x)=\iota=0. 
	\end{equation*}
holds. 
\end{lema}

\begin{proof}
	We first need to prove that 
		\begin{equation*}
			\lim_{x\to\infty}e^{-2nx}I_0(2nx)=0,
		\end{equation*}
	holds, for a fixed $n$ and $I_0$ being the modified Bessel function of the first kind of
	order $0$.
	The power series expansion for modified Bessel functions of the first kind of order $0$ is 
		\begin{equation*}
			I_0(x)=\sum^{\infty}_{k=0}{\frac{x^{2k}}{2^{2k}(k!)^2}},
		\end{equation*} 
	so for a fixed $n$ we have
		\begin{equation*}
			I_0(2nx)=\sum^{\infty}_{k=0}{\frac{(nx)^{2k}}{(k!)^2}}
		\end{equation*}
	and 
		\begin{equation*}
			e^{-2nx}\cdot I_0(2nx)=e^{-2nx}\cdot\sum^{\infty}_{k=0}{\frac{(nx)^{2k}}{(k!)^2}}=\varphi_n(x).
		\end{equation*}
	We now use Lebesgue's dominated convergence theorem and the integral expression
		\begin{align*}
			I_0(2nx)&=\frac{1}{\pi}\int^{1}_{-1}{e^{-2ntx}\cdot \frac{1}{\sqrt{1-t^2}}dt}\\
			e^{-2nx}\cdot I_0(2nx)&=\frac{1}{\pi}\int^{1}_{-1}{e^{-2nx(1+t)}\cdot \frac{1}{\sqrt{1-t^2}}dt},
		\end{align*}
	for $n$ fixed and we conclude that $\sigma_n(x)\to 0$, as $x\to \infty$, because we see from above that $e^{-2nx}\cdot I_0(2nx)\to 0$, for $x\to\infty$.  
\end{proof}

\begin{cor}\label{cor6.3}
The new Chebyshev-Gr\"{u}ss-type inequality for the  Mirakjan-Favard-Sz$\acute{a}$sz operator  is:
	\begin{equation*}
		\left|T(f,g;x)\right|\leq \frac{1}{2}\cdot osc_{M_n}(f)\cdot osc_{M_n}(g),
	\end{equation*}	
where $f,g\in C_b[0,\infty)$, $osc_{M_n}(f)=\sup\{\left|f_k-f_l\right|:0\leq k<l<\infty\}$ and a similar definition applying to $g$. 
\end{cor}

\section{Chebyshev-Gr\"{u}ss-type inequalities for Baskakov operators}

In the book of F. Altomare and M. Campiti \cite{altomare:1994} (Sect. 5.3.10), the classical positive, linear Baskakov operators $(A_n)_{n\in\N}$ are defined as follows:
	\begin{equation*}
		A_n(f;x):=\sum^{\infty}_{k=0}{f\left(\frac{k}{n}\right)\binom{n+k-1}{k}\frac{x^k}{(1+x)^{n+k}}},
	\end{equation*}
for every $f\in E_2$, $x\in [0,\infty)$ and $n\geq 1$. 

\subsection{A new Chebyshev-Gr\"{u}ss-type inequality for Baskakov operators}
The procedure in this subsection completely parallels that of Section 6.1.
We set 
	\begin{equation*}
		\vartheta_n(x):=\frac{1}{(1+x)^{2n}}\sum^{\infty}_{k=0}{\binom{n+k-1}{k}^2\left(\frac{x}{1+x}\right)^{2k}},\text{ for }x\geq 0.
	\end{equation*}
We need to find the infimum:
	\begin{equation*}	
	\inf_{x\geq 0}\vartheta_n(x):=\epsilon\geq 0.
	\end{equation*}
Because $\vartheta_n(x)\geq \epsilon$, we obtain the following result.

\begin{thm}
	For the Baskakov operator one has
	\begin{align*}
		\left|T(f,g;x)\right|&\leq \frac{1}{2}\left(1-\vartheta_n(x)\right)\cdot osc_{A_n}(f)\cdot osc_{A_n}(g)\\
		&\leq \frac{1}{2}(1-\epsilon)\cdot osc_{A_n}(f)\cdot osc_{A_n}(g),
	\end{align*}	
where $f,g\in C_b[0,\infty)$, $osc_{A_n}(f)=\sup\{\left|f_k-f_l\right|:0\leq k<l<\infty\}$, $f_k:=f\left(\frac{k}{n}\right)$ and a similar definition applying to $g$. 	
\end{thm}

\begin{lema} 
	The relation 
		$\inf_{x\geq 0}\vartheta_n(x)=\epsilon=0$
	holds for all $n\geq 1$. 
\end{lema}

\begin{proof}
	In \cite{berdysheva:2007} the following functions were defined. 
	For $I_c=[0,\infty) (c\in \R, c\geq 0)$, $n>0$, $k\in\N_0$ and $x\in I_c$, we have
		\begin{equation*}
			p^{[c]}_{n,k}(x):=(-1)^k\binom{-\frac{n}{c}}{k}(cx)^k(1+cx)^{-\frac{n}{c}-k}, c\neq 0. 
		\end{equation*}
	For $c=1$, we get
		\begin{equation*}
			p^{[1]}_{n,k}(x)=p_{n,k}(x)=(-1)^k\binom{-n}{k}x^k(1+x)^{-n-k}
			=\binom{n+k-1}{k}x^k(1+x)^{-n-k}=:a_{n,k}(x),
		\end{equation*}
	so we obtain the fundamental functions of the Baskakov operator.
	The following kernel function was defined in \cite{berdysheva:2007}:
		\begin{equation*}
			T_{n,c}(x,y)=\sum^{\infty}_{k=0}{p^{[c]}_{n,k}(x)\cdot p^{[c]}_{n,k}(y)}, \text{ for } x,y\in I_c.
		\end{equation*}
	We are interested in the case $c=1$ and $x=y$, so the above kernel becomes
		\begin{equation}
			T_{n,1}(x,x)=\sum^{\infty}_{k=0}{p^2_{n,k}(x)}=\sum^{\infty}_{k=0}{a^2_{n,k}(x)}=:\vartheta_n(x).
		\end{equation}
	For $n=1$, we get 
	\begin{equation*}		
		\vartheta_1(x)=T_{1,1}(x,x)=\frac{1}{(1+x)^2}\sum^{\infty}_{k=0}{\left(\frac{x}{1+x}\right)^{2k}}
		=\frac{1}{1+2x} \longrightarrow 0, \text{ for }x\to \infty.
	\end{equation*}
	For $n>1$, 
		\begin{equation*}
			T_{n,1}(x,x)=\frac{1}{\pi}\int^{1}_{0}{(\phi (x,x,t))^n\frac{dt}{\sqrt{t(1-t)}}},
		\end{equation*}
	where, for $\phi(x,x,t)=[1+4x(1-t)+4x^2(1-t)]^{-1}$, it holds:
		\begin{equation*}
			0<\phi(x,x,t)\leq 1, \forall t\in [0,1], \forall x\geq 0.
		\end{equation*}
	Therefore
		\begin{equation}\label{eq*}
			T_{2,1}(x,x)\geq T_{3,1}(x,x)\geq T_{4,1}(x,x)\geq \ldots \geq 0, \forall x\geq 0.
		\end{equation}
	Now for $n=2$, we have
		\begin{equation*}
			T_{2,1}(x,x)=\sum^{\infty}_{k=0}{p^2_{2,k}(x)}=\frac{1}{(1+x)^4}\sum^{\infty}_{k=0}{(k+1)^2\left(\frac{x}{1+x}\right)^{2k}}.
		\end{equation*}
	Let $\left(\frac{x}{1+x}\right)^2=y$. Then
		\begin{equation*}
			\sum^{\infty}_{k=0}{(k+1)^2 y^{k}}=\frac{1+y}{(1-y)^3}.
		\end{equation*}
	Thus 
		\begin{equation}\label{eq**}		
		T_{2,1}(x,x)=\frac{2x^2+2x+1}{(2x+1)^3}\to 0,\text{ for } x\to \infty. 
		\end{equation}
	For $n\geq 3$ it holds from (\ref{eq*}) that $0\leq T_{n,1}(x,x)\leq T_{2,1}(x,x)$. Combining this with (\ref{eq**}), we get
		\begin{equation*}
			\lim_{x\to 0}T_{n,1}(x,x)=0, \forall n\geq 1, 
		\end{equation*} 
	and so the proof is finished. 
\end{proof}

An inequality analogous to the one in Corollary \ref{cor6.3} is now immediate.

\section{Chebyshev-Gr\"{u}ss-type inequalities for Bleimann-Butzer-Hahn operators}

In the same book \cite{altomare:1994} (Sect. 5.2.8), the Bleimann-Butzer-Hahn operators are also presented. 
For every $n\in\N$ the positive linear operator $H_n:C_b([0,\infty))\to C_b([0,\infty))$ is defined by 
	\begin{equation*}
		H_n(f;x):=\frac{1}{(1+x)^n}\sum^{n}_{k=0}{f\left(\frac{k}{n-k+1}\right)\binom{n}{k}x^k},
	\end{equation*}
for every $f\in C_b([0,\infty))$, $x\geq 0$, $n\in\N$, $n\geq 1$.

\subsection{A new Chebyshev-Gr\"{u}ss-type inequality for Bleimann-Butzer-Hahn operators}

We set 
	\begin{equation*}
		\psi_n(t)=\frac{1}{(1+t)^{2n}}\sum^{n}_{k=0}{\binom{n}{k}^2t^{2k}},
	\end{equation*}
for  $t\geq 0$. We make a change of variable, namely set $x=\frac{t}{t+1}\in[0,1)$. Then we get
	\begin{align*}
		\psi_n(t)&=\sum^{n}_{k=0}{\binom{n}{k}^2\left(\frac{t}{t+1}\right)^{2k}\left(\frac{1}{t+1}\right)^{2n-2k}}\\
		&=\sum^{n}_{k=0}{\binom{n}{k}^2x^{2k}(1-x)^{2n-2k}}.
	\end{align*}
So $\psi_n(t)=\varphi_n(x)$, i.e., $\inf_{t\geq 0}\psi_n(t)=\inf_{x\in[0,1]}\varphi_n(x)=\frac{1}{4^n}\binom{2n}{n}$, as shown in Lemma \ref{lc3}.

This leads to

\begin{thm}
The new Chebyshev-Gr\"{u}ss-type inequality in this case is:
	\begin{equation}
		\left|T(f,g;x)\right|\leq \frac{1}{2}\left(1-\frac{1}{4^n}\binom{2n}{n}\right)\cdot osc_{H_n}(f)\cdot osc_{H_n}(g),
	\end{equation}
with 	$f,g\in C_b[0,\infty)$, $x\in [0,\infty)$ and  
	\begin{equation*}
		osc_{H_n}(f):=\sup \left\{\left|f_k-f_l\right|:0\leq k<l\leq n\right\},
	\end{equation*}
for $f_k:=f\left(\frac{k}{n-k+1}\right)$ and a similar definition applying to $g$. 
\end{thm}
\section{Chebyshev-Gr\"{u}ss-type inequalities for King-type operators}

P. P. Korovkin \cite{korovkin:1960} introduced in 1960 a theorem saying that if $\{L_n\}$ is a sequence of positive linear operators on $C[a,b]$, then 
	\begin{equation*}
		\lim_{n\to \infty}L_n(f)(x)=f(x)
	\end{equation*}
for each $f\in C[a,b]$ if and only if 
	\begin{equation*}
		\lim_{n\to \infty}L_n(e_i(x))=e_i(x)
	\end{equation*}
for the three functions $e_i(x)=x^{i},\ i=0,1,2$. There are a lot of well-known operators, like the Bernstein polynomials, the Mirakjan-Favard-Sz$\acute{a}$sz and the Baskakov operators, that preserve $e_0$ and $e_1$ (see \cite{king:2003}). However, these operators do not reproduce $e_2$. We are now interested in a non-trivial sequence of positive linear operators $\{L_n\}$ defined on $C[0,1]$, that preserve $e_0$ and
$e_2$:
	\begin{equation*}
		L_n(e_0)(x)=e_0(x) \text{ and } L_n(e_2)(x)=e_2(x),\ n=0,1,2,\ldots.
	\end{equation*}
In \cite{king:2003} J.P. King defined the King-type operator as follows. 

\begin{defin}(see \cite{king:2003})
	Let $\{r_n(x)\}$ be a sequence of continuous functions with $0\leq r_n(x)\leq 1$.
	Let $V_n:C[0,1]\to C[0,1]$ be defined by
		\begin{align*}
			V_n(f;x)&=\sum^{n}_{k=0}{\binom{n}{k}(r_n(x))^k(1-r_n(x))^{n-k}f\left(\frac{k}{n}\right)}\\
			&=\sum^{n}_{k=0}{v_{n,k}(x)\cdot f\left(\frac{k}{n}\right)},
		\end{align*}
	for $f\in C[0,1]$, $0\leq x\leq 1$. $v_{n,k}$ are the fundamental functions of the $V_n$ operator. 
\end{defin} 

\begin{rema}
For $r_n(x)=x,\ n=1,2,\ldots$, the positive linear operators $V_n$ given above reduce to the Bernstein
operators.
\end{rema}

\begin{prop}[Properties of $V_n$]\mbox{}
	\begin{enumerate}[i)]
		\item
			$V_n(e_0)=1$ and $V_n(e_1;x)=r_n(x)$;
		\item
			$V_n(e_2;x)=\frac{r_n(x)}{n}+\frac{n-1}{n}(r_n(x))^2$;
		\item
			$\lim_{n\to\infty}V_n(f;x)=f(x)$ for each $f\in C[0,1]$, $x\in [0,1]$, if and 
			only if  
				\begin{equation*}
				\lim_{n\to \infty}r_n(x)=x.
				\end{equation*}
	\end{enumerate}
\end{prop}

For special ("right") choices of $r_n(x)=r^*_n(x)$, J. P. King showed in \cite{king:2003} that the following theorem holds.
\begin{thm} (see Theorem 1.3. in \cite{gopi:2005})
	Let $\{V^*_n\}_{n\in\N}$ be the sequence of operators defined before with
	\begin{equation*}
	r^*_n(x):=
	\begin{cases}
		r^*_1(x)=x^2 & \text{, for } n=1, \\
		r^*_n(x)=-\frac{1}{2(n-1)}+\sqrt{\left(\frac{n}{n-1}\right)x^2+\frac{1}{4(n-1)^2}} & \text{, for } n=2,3,\ldots.
	\end{cases}
	\end{equation*}
	Then we get $V^*_n(e_2;x)=x^2$, for $n\in\N$, $x\in[0,1]$ and $V^*_n(e_1;x)\neq e_1(x)$. 	$V^*_n$ is not a polynomial operator.  
\end{thm}

The fundamental functions of this operator, namely
	\begin{equation*}
		v^*_{n,k}(x)=\binom{n}{k}(r^*_n(x))^k(1-r^*_n(x))^{n-k}
	\end{equation*}
satisfy $\sum^{n}_{k=0}{v^*_{n,k}(x)}=1$, for $n=1,2,\ldots$. 

\begin{prop}[Properties of $r^*_n$]\mbox{}
	\begin{enumerate}[i)]
		\item
			$0\leq r^*_n(x)\leq 1$, for $n=1,2,\ldots$, and $0\leq x \leq 1$.
		\item
			$\lim_{n\to\infty}r^*_n(x)=x$ for $0\leq x\leq 1$.
	\end{enumerate}
\end{prop}

\subsection{The classical Chebyshev-Gr\"{u}ss-type inequality for King-type operators}

The second moments of the special King-type operators $V^{*}_{n}$ are given by
	\begin{equation*}
		V^*_n((e_1-x)^2;x)=2x(x-r^*_n(x)),
	\end{equation*}
so we discriminate between two cases. 

The first case is $n=1$, so $r^*_n(x)=x^2$ and the second moment is
	\begin{equation*}
		V^*_1((e_1-x)^2;x)=2x^2(1-x),
	\end{equation*}

so the classical Chebyshev-Gr\"{u}ss-type inequality is given as follows. 

\begin{thm}
	For $L=V^*_1$, we have the inequality:
		\begin{equation*}
	\left|T(f,g;x)\right|\leq \frac{1}{4}\widetilde{\omega}\left(f;2x\sqrt{2(1-x)}\right)\cdot \widetilde{\omega}\left(g;2x\sqrt{2(1-x)}\right).
	\end{equation*}
\end{thm}

For the second case, $n=2,3,\ldots$, we have
	\begin{equation*}
		r^*_n(x)=-\frac{1}{2(n-1)}+\sqrt{\left(\frac{n}{n-1}\right)x^2+\frac{1}{4(n-1)^2}},
	\end{equation*}
so the second moments look more complicated:
	\begin{equation*}
		V^*_n((e_1-x)^2;x)=2x(x-r^*_n(x)).
	\end{equation*}

In this case we get the following:

\begin{thm}
	For $L=V^*_n(x)$ and $n=2,3,\ldots$, the inequality
		\begin{equation*}
		\left|T(f,g;x)\right|\leq \frac{1}{4}\widetilde{\omega}\left(f;2\sqrt{2x(x-r^*_n(x))}\right)\cdot \widetilde{\omega}\left(g;2\sqrt{2x(x-r^*_n(x)}\right)
		\end{equation*}
	holds.
\end{thm}

\subsection{A new Chebyshev-Gr\"{u}ss-type inequality for King-type operators}

We need $\sum^{n}_{k=0}{(v^*_{n,k}(x))^2}$ to be minimal. Let $\varphi_{n}(x):=\sum^{n}_{k=0}{(v^*_{n,k}(x))^2}$. 

For $n=1$, we have that 
	\begin{align*}	
	\varphi_{1}(x)&=\sum^{1}_{k=0}{(v^*_{1,k}(x))^2}\\
					&=(v^*_{1,0}(x))^2+(v^*_{1,1}(x))^2\\
					&=2x^4-2x^2+1
	\end{align*}
and this attains its minimum for $x=\frac{\sqrt{2}}{2}$. This minimum is 
	\begin{equation*}
		\varphi_{1}\left(\frac{\sqrt{2}}{2}\right)=\frac{1}{2}.
	\end{equation*}

\begin{thm}
The new Chebyshev-Gr\"{u}ss-type inequality for $n=1$ then looks as follows:

	\begin{align*}
		\left|T(f,g;x)\right|&\leq \frac{1}{4}\cdot osc_{V^*_1}(f)\cdot osc_{V^*_1}(g)\\
									  &=\frac{1}{4}\cdot \left|f_0-f_1\right|\cdot\left|g_0-g_1\right|.
	\end{align*}
\end{thm}
	
For $n=2,3,\ldots$, the problem of finding the minimum is more difficult, since

\begin{align*}
	\varphi_{n}(x)&=\sum^{n}_{k=0}{(v^*_{n,k}(x))^2}\\
	&=\sum^{n}_{k=0}{\binom{n}{k}^2(r^*_{n}(x))^{2k}(1-r^*_{n}(x))^{2(n-k)}}.
\end{align*}

In any case, the estimate 
	\begin{equation*}
		\varphi_{n}(x)=\sum^{n}_{k=0}{(v^*_{n,k}(x))^2}\geq \frac{1}{n+1}
	\end{equation*}
holds, for $x\in [0,1]$ and $n=2,3,\ldots$. As a proof for this,
	\begin{equation*}
		\sqrt{\frac{\sum^{n}_{k=0}{v^*_{n,k}(x))^2}}{n+1}}\geq \frac{\sum^{n}_{k=0}{v^*_{n,k}(x)}}{n+1}=\frac{1}{n+1}.
	\end{equation*}
Then we get
	\begin{equation*}
		1-\sum^{n}_{k=0}{(v^*_{n,k}(x))^2}\leq 1-\frac{1}{n+1}=\frac{n}{n+1}.
	\end{equation*}

\begin{thm}
For $n=2,3,\ldots$ there holds 
	\begin{equation*}
		\left|V^*_n(fg)(x)-V^*_n(f;x)\cdot V^*_n(g;x)\right|\leq \frac{n}{2(n+1)}\cdot osc_{V^*_n}(f)\cdot osc_{V^*_n}(g).
	\end{equation*}
\end{thm}

\bigskip
 \noindent
 $
 \begin{array}{l}
 \textrm{Heiner Gonska}\\
 \textrm{University of Duisburg-Essen} \\
 \textrm{Faculty of Mathematics } \\
 \textrm{Forsthausweg 2} \\
 \textrm{47057 Duisburg} \\
 \textrm{Germany}\\
\textrm{e-mail: heiner.gonska@uni-due.de} 
\end{array} 
 $
$
\begin{array}{lll}
 \textrm{Ioan Ra\c{s}a}\\
 \textrm{Technical University of Cluj-Napoca} \\
 \textrm{Department of Mathematics } \\
 \textrm{Str. C. Daicoviciu, 15} \\
 \textrm{RO-400020 Cluj-Napoca} \\
 \textrm{Romania}\\
\textrm{e-mail: Ioan.Rasa@math.utcluj.ro} 
\end{array} 
$
\\
 $
 \begin{array}{ll}
 \textrm{Maria-Daniela Rusu}\\
 \textrm{University of Duisburg-Essen} \\
 \textrm{Faculty of Mathematics } \\
 \textrm{Forsthausweg 2} \\
 \textrm{47057 Duisburg} \\
 \textrm{Germany}\\
\textrm{e-mail: maria.rusu@uni-due.de} 
\end{array} 
$

\end{document}